\documentclass[12pt]{article}

\pdfoutput=1

\usepackage{amsmath,amssymb}
\usepackage{graphicx}
\newcommand{\dd}[2]{\frac{\textstyle\partial #1}{\textstyle\partial #2}}
\renewcommand{\dot}[1]{\overset{\text{\large .}}{#1}}
\newcommand{\textfrac}[2]{\ensuremath{\frac{\textstyle #1}{\textstyle #2}}}
\begin{document}
\thispagestyle{empty}
\begin{center}
\large\bf The Equations of Motion of a Charged Particle
in the Five-Dimensional Model of the General Relativity Theory
with the Four-Dimensional Nonholonomic Velocity Space
\end{center}
\begin{center}
\textit{\large Victor R. Krym\footnote{E-mail: vkrym2007@rambler.ru}, 
Nickolai N. Petrov}\\
\end{center}
\begin{center}
St.-Petersburg State University, Department of Mathematics and Mechanics,\\
St.-Petersburg, 198504 Russia
\end{center}

\begin{abstract}
We consider the four-dimensional nonholonomic distribution defined by the 
4-potential of the electromagnetic field on the manifold.
This distribution has a metric tensor with the Lorentzian
signature $(+,-,-,-)$, therefore, the causal structure appears 
as in the general relativity theory.
By means of the Pontryagin's maximum principle
we proved that the equations of the horizontal geodesics for this
distribution are the same as the equations of motion of a charged
particle in the general relativity theory. 
This is a Kaluza -- Klein problem of classical and quantum physics
solved by methods of sub-Lorentzian geometry.
We study the geodesics sphere which appears in a constant magnetic field
and its singular points.
Sufficiently long geodesics are not optimal solutions of the variational problem and
define the nonholonomic wavefront.
This wavefront is limited by a convex elliptic cone.
We also study variational principle approach to the problem.
The Euler -- Lagrange equations are the same as those obtained by the Pontryagin's
maximum principle if the restriction of the metric tensor on the distribution is the same.
\end{abstract}


AMS Subj. Class. (MSC): 37J60, 37J55, 53B30, 53B50, 53D50, 58A30, 70F25, 81S05


Journal-ref: Vestnik Sankt-Peterburgskogo Universiteta, Ser. 1. 
Matematika, Mekhanika, Astronomiya, 2007, N1, pp. 62--70.

\section{Nonholonomic model in sub-Lorentzian geometry} 

Distribution on a smooth manifold $M$ is a family of subspaces 
$\mathcal A(x)\subset T_x M$ (of the tangent bundle of the manifold) with the 
same dimension, smoothly parametrized by the points of manifold
 \cite{Gromov,VerGer,Dobr,Newm}. 
For each $x\in M$ on $\mathcal A(x)$ the bilinear form 
$\langle u, u\rangle_x$ is defined which can be equally described by the metric 
tensor  $g_{ij}(x)$.
The metric tensor of a distribution smoothly depends of the point of the 
manifold.
The absolutely continous path $x(t)$ is called {\it horizontal}, iff
$x'(t)\in\mathcal A(x(t))$ for almost any $t$.
In sub-Riemannian geometry the metric tensor $g_{ij}(x)$ of the distribution 
$\mathcal A$ 
is positively defined. In this geometry we study the {\it shortest} horizontal 
curves connecting two given points. 
If the distribution is completely nonholonomic and the Riemannian manifold is 
full,
then any two points of the manifold can be connected by the shortest horizontal 
geodesic.
On manifolds and also on distributions with metric tensor the geodesic is the 
Euler -- Lagrange solution for the length functional
$J(x(\cdot),u(\cdot))= \int_{t_0}^T \langle u(t),u(t)\rangle^{1/2}_{x(t)}
\, dt$.

In sub-Lorentzian geometry the metric tensor $g_{ij}(x)$ of the distribution 
$\mathcal A$
has the signature  $(+,-,\ldots\!,-)$. In this geometry we study the {\it 
longest} horizontal curves connecting two given points \cite{Beem}.
The equations of motion of a charged particle in the electromagnetic and 
gravitational
fields can be obtained as the solutions of the variational problem for the
distribution  \cite{Kr}.
Consider the 5-dimensional manifold $M^5$ with the 4-dimensional distribution
$\mathcal A$ defined by means of the differential form
$\omega(x)=\sum\limits_{i=0}^3 A_i(x) dx^i + dx^4$, 
where the covector  $(A_i)_{i=0,\ldots\!,3}$ is the 4-potential of the 
electromagnetic field.
Consider the length functional
\begin{equation}
J(x(\cdot),u(\cdot))= \int_{t_0}^T L(x(t),u(t)) \, dt,
\end{equation}
where $L(x,u)=m\langle u, u\rangle^{1/2}_x$ is the pseudonorm of the vector $u$,
$m$ is the mass of the particle.
Consider the problem of maximizing the length functional on the set of
horizontal curves connecting two given points.
Solution of this problem as we proove in this paper leads to equations
\begin{equation}\label{genreleqn}
\frac{d}{dt} \dd{L}{u^k} - \dd{L}{x^k} + p_4 \sum_{j=0}^3 F_{jk} u^j = 0,
\quad k=0,\ldots\!,3,
\end{equation}
where $F_{jk}$ are defined by  \eqref{tensorF}.
These equations are the same as the equations of motion of a particle
with the charge  $p_4$ in the electromagnetic and gravitational fields of the 
general relativity theory.
This wording of the problem is different from the classical problem of physics 
\cite{Landau,Gray} of minimizing the action functional on a 4-dimensional manifold:
\begin{equation}
S(x(\cdot),u(\cdot))=
-m c \int_0^{t_1} \langle u(t),u(t)\rangle^{1/2}_{x(t)} \, dt - 
\frac{e}{c} \int_0^{t_1} \sum_{i=0}^3 A_i(x(t)) u^i(t) \, dt,
\end{equation}
where $e$ is the charge of the particle, $c$ is the speed of light.
This paper is devoted to the problem of unification of electromagnetic,
gravitational and other fundamental interactions \cite{Rumer,BailinL,Sriv}.

In this model the particles are considered on the 5-dimensional manifold
$M^5$ with the special smooth structure.
Allowed are only smooth coordinate transformations for which
\begin{equation} \label{coordinate}
\dd{y^i}{x^4} = 0, \quad i=0,\ldots\!,3, \quad
\dd{y^4}{x^4} = 1.
\end{equation}
Components $\dd{y^4}{x^i}$, $i=0,\ldots\!,3$, can be interpreted as gauge 
transformations \eqref{gauge}.
Since the coordinate transformations are smooth,
$\dd{}{x^4}\dd{y^i}{x^j}=\dd{}{x^j}\dd{y^i}{x^4}=0$,
$i,j=0,\ldots\!,4$. 
Therefore the cylindricity condition
$\dd{v}{x^4}=0$ for any tensor field
$v$  is invariant.
In the tangent bundle $TM$ the subspace 
$\text{Lin}\,\{\partial_4\}$ is invariant.
In the cotangent bundle $T^*M$ the subspace 
$\text{Lin}\,\{dx^0,dx^1,dx^2,dx^3\}$ is invariant.
Assume that the metric tensor and the 4-potential of the electromagnetic field
do not depend of  $x^4$. 
This condition is invariant at the coordinate transformations with the property
\eqref{coordinate}.
The 4-dimensional distribution  $\mathcal A$ is the velocity space,
i.e. the set of all possible velocity vectors of the particle.
The metric tensor of the distribution
$\mathcal A$ defines 4-dimensional convex elliptic cone,
completely similar to the cone of future of the general relativity theory. 
Causal structure in the proposed theory is similar to the causal structure
of the general relativity theory \cite{KrM,KrP}. 
In the presence of the electromagnetic field the manifold $\mathcal A$ 
is nonholonomic.

The electromagnetic field is defined in physics by the antisymmetric tensor 
$(F_{jk})_{j,k=0,\ldots\!,3}$
\begin{equation} \label{tensorF}
(F_{jk})_{j,k=0,\ldots\!,3}= 
\begin{pmatrix}
0 & E_x & E_y & E_z \\
-E_x & 0 & -H_z & H_y \\
-E_y & H_z & 0 & -H_x \\
-E_z & -H_y & H_x & 0
\end{pmatrix}.
\end{equation}
Here $E$ is the 3-dimensional vector of the electric field tension,
$H$ is the 3-dimensional vector of the magnetic field tension. 
There is the 4-potential $(A_i)_{i=0,\ldots\!,3}$ which defines the tensor of 
the electromagnetic field:
\begin{equation}
F_{jk}=\dd{A_k}{x^j}-\dd{A_j}{x^k} \ , \quad j,k=0,\ldots\!,3.
\end{equation}
4-Potential is not defined unambiguously.
It is defined with an arbitrary gauge transformation
$A_i\mapsto A_i+\dd{f}{x^i}$, where $f$ is any smooth function.

Let $U$ be some coordinate neighbourhood on the smooth manifold $M^5$.
The 4-dimensional distribution $\mathcal A$ on $U$ can be defined by the 
differential form 
\begin{equation}
\omega_x=\sum_{i=0}^3 A_i(x) dx^i + dx^4.
\end{equation}
This definition is correct since at the coordinate transformation with the 
property
\eqref{coordinate} we have
\begin{equation}
\sum_{i=0}^3 A_i dx^i + dx^4 =
\sum_{j=0}^3 \left( \sum_{i=0}^3 A_i \dd{x^i}{y^j} + \dd{x^4}{y^j}
\right) dy^j + dy^4 =
\sum_{j=0}^3 \widetilde A_j dy^j + dy^4.
\end{equation}
Therefore the coordinate transformations on the 5-manifold lead to
transformations of the 4-potential
\begin{equation}\label{gauge}
\widetilde A_j = \sum\limits_{i=0}^3 A_i \dd{x^i}{y^j} + \dd{x^4}{y^j},
\end{equation}
that include 4-dimensional coordinate transformations of the general relativity 
theory
and gauge transformations of the 4-potential of the electromagnetic field 
\cite{Kr} as in
 \cite{Landau}. 
Locally the distribution  $\mathcal A$ can be defined by the basis vector fields
$e_i=\dd{}{x^i}-A_i\dd{}{x^4}$, $i=0,\ldots\!,3$.
Vector fields on a manifold can be considered as operators of differentiation.
During the coordinate transformation with the property 
\eqref{coordinate} we have
\begin{equation}\label{coordbasis}
\widetilde e_j=\sum\limits_{i=0}^3 \dd{x^i}{y^j} e_i,
\end{equation}
where
$\widetilde e_j=\dd{}{y^j}-\widetilde A_j\dd{}{y^4}$, $j=0,\ldots\!,3$.
Therefore the 4-dimensional basis of the distribution
$\mathcal A$ is transformed \eqref{coordbasis} exactly the same way as  
the coordinate vector fields $\left(\dd{}{x^i}\right)_{i=0,\ldots\!,3}$ 
in the general relativity theory.
The commutators of vector fields
$e_i$, $i=0,\ldots\!,3$, generate the tensor of the electromagnetic field: 
$[e_i, e_j]=-F_{ij}\dd{}{x^4}$.
If  $F_{ij}\ne 0$, then the distribution $\mathcal A$ is completely 
nonholonomic.
Allowed curves on the distribution are absolutely continuous and
satisfy almost everywhere the horizontality condition
\begin{equation} \label{horizontal}
\omega_{\gamma(t)}(\gamma'(t))=0.
\end{equation}

In this paper we consider the optimization problem for the length functional\footnote{
Note that nonholonomic variational problem is different from the 
nonholonomic mechanical problem.}
on the 5-dimensional manifold with the 4-dimensional distribution
$\mathcal A$. We prove that the solutions of the variational problem for this
distribution satisfy the equations of motion for the charged particle of the
general relativity theory. We consider also the geodesic sphere
for the distribution $\mathcal A$.
{\it The geodesic sphere} with the radius $r$ and center  $x_0$ is the set of 
points at distance  $r$ from the point $x_0$.
To obtain the geodesic sphere we consider the set of solutions
of the variational problem starting at point $x_0$.
We constructed the 4-dimensional projection of the geodesic sphere
for the particle moving in the constant magnetic field.

\section{The equations of motion of a charged particle}
\label{proof}

The optimization problem for the length functional on the set of the
horizontal curves connecting two given points is one of the main problems
of the theory of optimal control \cite{Vasil,KrotG}.
The velocity vector $x'(t)$ of horizontal curves belong to the distribution 
$\mathcal A$,
hence we consider the metric tensor of the distribution $\mathcal A$. 
If the vector fields $e_i(x)$, $i=0,\ldots\!,3$ form the basis of the 
distribution 
$\mathcal A(x)$, then its metric tensor is
$g_{ij}(x)=\langle e_i(x), e_j(x) \rangle_x$, $i,j=0,\ldots\!,3$.
Consider the length functional
\begin{equation} \label{length}
J(x(\cdot),u(\cdot))= \int_{t_0}^T L(x(t),u(t)) \, dt,
\end{equation}
where
\begin{equation}
L(x,u)= m \left( \sum\limits_{i,j=0}^3 g_{ij}(x) \, u^i \, u^j \right)^{1/2}
\end{equation}
and $m$ is the mass of the particle. 
Since the metric tensor has the signature
$(+,-,-,-)$, the velocity vector
$\textfrac{dx}{dt}$ belong to the cone of future  
$V(x) = \{ v\in\mathcal A(x) \mid \langle v, v \rangle_x >0,\  v^0>0 \}$.
This set depends of the point $x$. 
Hence let us consider another problem.

All spaces $\mathcal A(x)$ are linearly isomorphic. Therefore we can consider 
the cone 
$\overline V(x) = \{ u\in\mathbb R^4 \mid 
\sum\limits_{i,j=0}^3 g_{ij}(x) u^i u^j >0,\  u^0>0 \}$.
For some neighbourhood of the point $x_0$
the inersection of all sets  $\overline V(x)$ is nonempty and has a nonempty 
interior $Q$. In this neighbourhood of the point  $x_0$ consider another problem
for which $u(t)\in Q$.
If $(x(t),u(t))$ is an optimal process for the problem \eqref{length},
for which $u(t)\in\overline V(x(t))$, then this pair is locally an optimal 
process in 
the case  $u(t)\in Q$.
From the horizontality condition  \eqref{horizontal} follows that 
for the velocity vector of the particle 
\begin{equation} \label{fxu}
\left\{ \begin{aligned}
\frac{dx^k}{dt} &= u^k, \quad k=0,\ldots\!,3 \\
\frac{dx^4}{dt} &= -\sum\limits_{j=0}^3 A_j(x) \, u^j.
\end{aligned} \right.
\end{equation}
Let us consider the maximization problem for the length functional on the set of 
horizontal curves connecting two given points
$x^i(t_0)=x^i_0$ and
$x^i(T)=x^i_T$,  $i=0,\ldots\!,4$.
The Hamilton -- Pontrjagin function for the functional \eqref{length},
\eqref{fxu} has the form
\begin{equation}
H(x,u,p,a_0) = -a_0 L(x,u) + 
\sum_{j=0}^3 p_j u^j - p_4 \sum_{j=0}^3 A_j(x) u^j.
\end{equation}
Due to Pontrjagin's maximum principle \cite{Pont, Phil}
there is  $a_0\geq 0$ and vector-function
$p_k(t)$, $k=0,\ldots\!,4$, $t_0\leq t \leq T$,
satisfying the conjugate system of linear differential equations
\begin{equation} \label{dpdt}
\frac{dp_k}{dt} = -\left.\dd{H}{x^k}\right|_{(x(t),u(t),p(t),a_0)}
= a_0 \dd{L}{x^k} 
+ p_4 \sum_{j=0}^3 \dd{A_j}{x^k} u^j(t),
\quad k=0,\ldots\!,4,
\end{equation}
where $(x(t),u(t))$ is the optimal process for the problem  \eqref{length}.
Since the metric tensor and 4-potential of the electromagnetic field
do not depend of the coordinate $x^4$, the appropriate momentum component
is the integral of motion (is conserved): $p_4(t)=\text{const}$.
From the equation \eqref{genreleqn} follows that $p_4$ is the length of the 
particle.
For almost any $t\in [t_0,T]$ the function $u\mapsto H(x(t),u,p(t),a_0)$
reaches its upper bound on the set  $Q$ at $u=u(t)$.
If this maximum is reached at an inner point, then
$\dd{H}{u^k}=0$, i.e.
\begin{equation} \label{Hmax}
p_k - p_4 A_k = a_0 \dd{L}{u^k}, \quad k=0,\ldots\!,3.
\end{equation}
Since the length functional does not depend on the parametrization of the path
$x(\cdot)$, we can assume that the optimal control belong to the pseudosphere
$L(x,u)=1$. The normal covector for this pseudosphere has the coordinates
$\dd{L}{u^k}$, $k=0,\ldots\!,3$.
The projection of the impulse vector $p_k$ on the subset $(u^0,\ldots\!,u^3)$
has the coordinates $p_k-p_4 A_k$. The Hamilton function reaches its maximum 
when 
the projection of the vector  $p_k$ and the normal for the unit pseudosphere are 
collinear.
The equation \eqref{Hmax} is exactly the condition for these two vector to be 
collinear.
To exclude the vector $p_k$ from the equations \eqref{dpdt} and \eqref{Hmax}
we use the identity
\begin{equation}
A_k(x(t))-A_k(x(t_0))=\int_{t_0}^t \sum_{j=0}^3 \dd{A_k}{x^j} u^j \, dt' , 
\quad k=0,\ldots\!,3.
\end{equation}
Then
\begin{equation}
a_0 \left( \dd{L}{u^k} - \int_{t_0}^t \dd{L}{x^k} \, dt' \right)
+ p_4 \int_{t_0}^t \sum_{j=0}^3 F_{jk} u^j \, dt' +
p_4 A_k(x(t_0))-p_k(t_0) =0.
\end{equation}
Parameters $a_0$ and $p_k$, $k=0,\ldots\!,4$,
can be multiplied by any positive constant.
With $a_0=1$ we obtain the equations of horizontal geodesics on the distribution 
$\mathcal A$. If we additionnaly assume continuity of functions $u(t)$, then
\begin{equation} \label{genrel}
\frac{d}{dt} \dd{L}{u^k} - \dd{L}{x^k} + p_4 \sum_{j=0}^3 F_{jk} u^j = 0,
\quad k=0,\ldots\!,3.
\end{equation}
These equations are identical with the equations of motion of a particle
with the charge  $p_4$ of the general relativity theory.

For $a_0=0$ we obtain the equations of abnormal geodesics.
Then at least one of the parameters $p_k$, $k=0,\ldots\!,4$, should be non-zero.
Since in this case \eqref{Hmax}
$p_k = p_4 A_k$, $k=0,\ldots\!,3$, then $p_4\ne 0$.
The solutions of the variational problem for $a_0=0$ has the form
\begin{equation} \label{singular}
\sum_{j=0}^3 F_{jk} u^j =0, \quad k=0,\ldots\!,3.
\end{equation}
The solutions of the variational problem
\eqref{length} -- \eqref{fxu} which do not satisfy Euler -- Lagrange equations
are called abnormal geodesics \cite{Petr,Bonn}.

Montgomery constructed a distribution in 3-dimensional space 
for which the commutator of the basis vector fields is
$2x^2\dd{}{x^3}$ \cite{Mont1,Mont2}.
In his example the abnormal geodesics are streight lines parallel with
the  $x^2$ axis.
In our model the equations of abnormal geodesics on the distribution
$\mathcal A$ has the form \eqref{singular}.
If $\det F=0$, then this equation can posess a non-trivial solution.

\section{The geodesic sphere for the particle in the magnetic field}

The geodesic sphere has a compicated structure \cite{VerGer,Kup,KupO}.
In the following part of the paper we assume that 
$M$ is the linear space $\mathbb R^5$ and the metric tensor of the distribution
$\mathcal A$ is diagonal,
$g = \text{Diag}(1,-1,-1,-1)$. 
The equations of motion of a charged particle in the electromagnetic field
\eqref{genrel} in the space with the diagonal metric tensor has the form
\begin{equation}
m\frac{du^k}{dt}-q\sum_{j=0}^3 F^k_{\ j\ } u^j =0, \quad k=0,\ldots\!,3,
\end{equation}
where  $m$ is the mass of the particle, $q$ is the charge of the particle, $F$ 
is the tensor of the electromagnetic field, $u^k$ are components of the 4-
velocity vector of the particle, 
$\textfrac{dx^k}{dt}=u^k$, $k=0,\ldots\!,3$.
The fiveth component of the velocity vector $\textfrac{dx^4}{dt}$ is determined 
by the horizontality condition.

Let us consider the motion of the particle in the magnetic field.
Assume that the tension of the magnetic field 
$H_x\ne 0$, and the tension of the electric field is zero.
Then $F_{23}=-F_{32}\ne 0$, and other components of the tensor of the 
electromagnetic field  $F$ are zero. In this coordinates the equations of motion 
and the horizontality condition \eqref{horizontal} has the form
\begin{equation} \label{eqnmove}
\left\{
\begin{aligned}
m\frac{du^2}{dt}-q H_x u^3 =0 \\
m\frac{du^3}{dt}+q H_x u^2 =0 \\
\frac{dx^4}{dt} + \sum_{k=0}^3 A_k u^k = 0.
\end{aligned}
\right.
\end{equation}
Other components of the 4-velocity vector are constant. Assume that 
$H_x=\text{const}$. Then from the equations \eqref{eqnmove} follows
that in the magnetic field the charged particle moves along the helicoid line.

Consider the projection of the geodesic sphere on the coordinates  
$(x^2,x^3,x^4)$.
In the plane  $(x^2,x^3)$ the particle moves along a circle.
Choose the following components of the 4-potential:
$A_0=0$, $A_1=0$, $A_2=\varphi x^3$, $A_3=0$.
Then $F_{23}=-\varphi$, $H_x=\varphi$. 
Designate $p=\textfrac{q}{m}\varphi$. 
Then the equations of motion of a particle in the subspace
$(x^2,x^3,x^4)$ has the form
\begin{equation}
\left\{
\begin{aligned}
&\frac{du^2}{dt}-p u^3 =0 \\
&\frac{du^3}{dt}+p u^2 =0 \\
&\frac{dx^4}{dt} + \varphi x^3 u^2 =0.
\end{aligned}
\right.
\end{equation}
\begin{figure}[t]
\includegraphics[width=7cm,height=7cm]{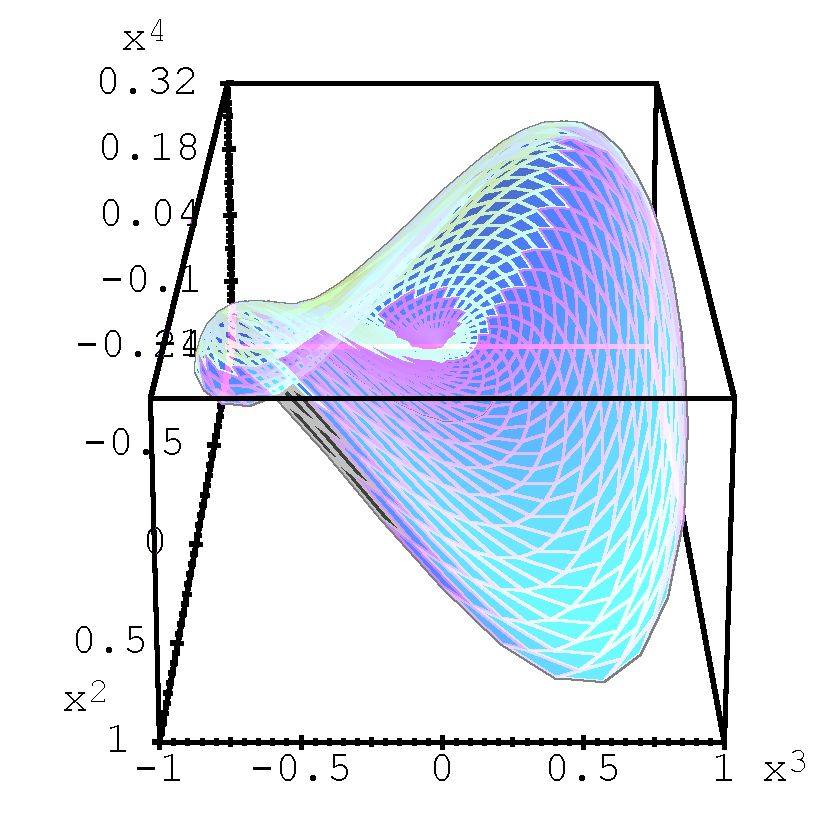}%
\includegraphics[width=7cm,height=7cm]{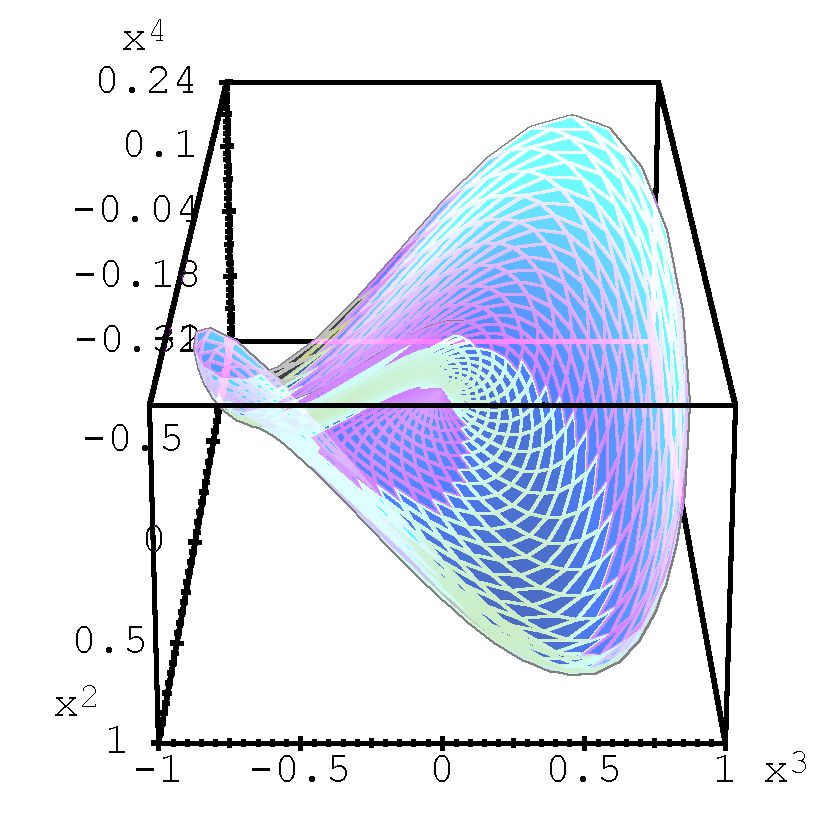}

\caption{Projection of the geodesic sphere on the coordinates $(x^2,x^3,x^4)$.}
\end{figure}
The norm of the velocity vector along the geodesic is constant.
Assume that $(u^2)^2+(u^3)^2=1$.
Then
\begin{equation} \label{x2x3}
\left\{
\begin{aligned}
x^2(t) &=-\frac{1}{p} \cos(\alpha+pt) + b_2\\
x^3(t) &=\frac{1}{p} \sin(\alpha+pt) + b_3\\
x^4(t) &= \frac{\varphi}{4p^2} (-2pt+\sin(2\alpha+2pt)) + 
\frac{\varphi}{p} b_3 \cos(\alpha+pt) + b_4.
\end{aligned}
\right.
\end{equation}
In the plane $(x^2,x^3)$ the charged particle is moving along the circle with 
the center 
$(b_2,b_3)$ and the radius $\textfrac{1}{|p|}$.
Let $b_2=(\cos\alpha)/p$ and $b_3=-(\sin\alpha)/p$.
Choose the constant $b_4$ so that $x^4(0)=0$.
Then
\begin{equation}\label{x2x3x4}
\left\{
\begin{aligned}
x^2(t) &=\frac{1}{p} (\cos\alpha - \cos(\alpha+pt))\\
x^3(t) &=\frac{1}{p} (\sin(\alpha+pt) - \sin\alpha)\\
x^4(t) &= \frac{\varphi}{4p^2} (-2pt + \sin(2\alpha+2pt) - \sin 2\alpha) +
\frac{\varphi}{p^2} \sin\alpha (\cos\alpha - \cos(\alpha+pt)).
\end{aligned}
\right.
\end{equation}
The length of the geodesic $\gamma:[0,t]\to\mathbb R^5$ is equal to $t$.

The equations \eqref{x2x3x4} define mapping
$(t,\alpha,p)\mapsto\gamma_{\alpha,p}(t)$.
In the plane $(x^2,x^3)$ the point makes one turn around the center 
at $|pt|=2\pi$. More long geodesics
($|pt|>2\pi$) are not the shortest paths.
Here one should consider the shortest paths, since the restriction of the metric 
tensor on the tangent bundle for the coordinate plane $(x^2,x^3,x^4)$
is positively defined.
The geodesic sphere of the radius  $s$ can be defined by the following formula:
\begin{equation}
B(s)=\{ \gamma_{\alpha,p}(t) \,\bigl|\,
\alpha\in [-\pi, \pi],\,
t=s,\, |pt|\leq 2\pi
\}.
\end{equation}
Fig.~1 shows the projection of the geodesic sphere of the radius
$s=1$ for the coordinates $(x^2, x^3, x^4)$ at $\varphi=1$, with
$p \in [-2\pi, 0]$ on the left and  $p \in [0, 2\pi]$ on the right.
If $|pt|\leq 2\pi$, the geodesics are the shortest paths.
The nonholonomic geodesic sphere of any radius has a pair of points such that
if we continue the geodesic after this point we enter the geodesic sphere with 
some lesser radius.
If $|pt|\geq 2\pi$, then the geodesic of length $t$ ceases to be optimal 
path.
All possible endpoints $\gamma_{\alpha,p}(t)$ of the geodesics of length
$t$ reach the $x^4$ axis at $pt=2\pi k$, $k\in\mathbb Z$,
this point being the intersection point for the considered surface.
The density of these points increases when we approach the origin.
This phenomena is called the nonholonomic wave front (fig.~2).
\begin{figure}[t]
\includegraphics[width=7cm,height=7cm]{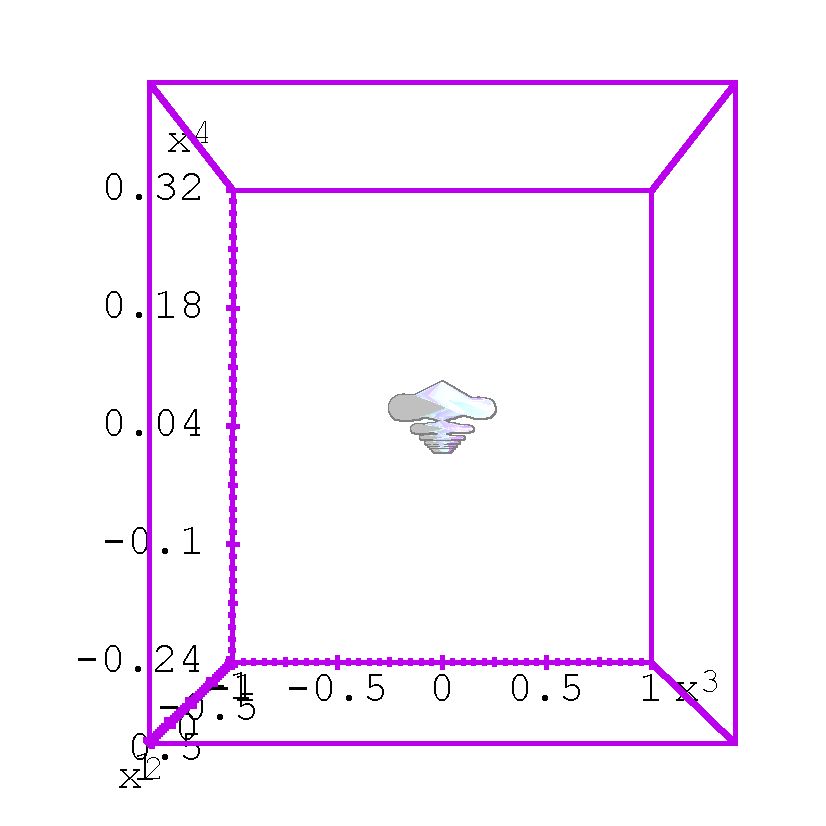}
\includegraphics[width=7cm,height=7cm]{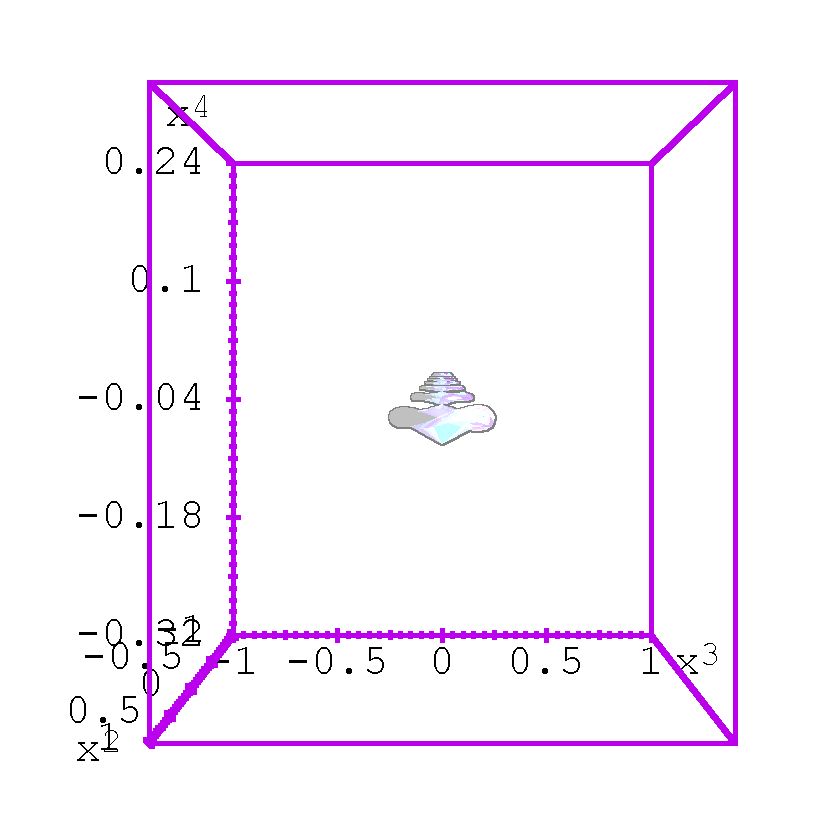}

\caption{The non-holonomic wave front defined by non-optimal geodesics.}
\end{figure}
Fig.~2 shows the projection of the surface defined by the geodesics of length 
$s=1$ for the coordinates $(x^2, x^3, x^4)$ at $\varphi=1$ with
$p \in [-8\pi,-\pi]$ on the left and $p \in [\pi,8\pi]$ on the right.
These geodesics cease to be optimal paths for 
$|p|>2\pi$.
For large enough $p$ this surface belongs to some cone with the apex 
at the origin. This cone has the axis of symmetry  $x^4$.
Indeed for  $p$ large enough
from the equations \eqref{x2x3x4} follows that
\begin{equation} \label{cone}
x^4(t)=-\frac{\varphi}{2p} s + 
O(\frac{1}{p^2}) \quad (p\to \infty).
\end{equation}
In the plane $(x^2,x^3)$ the particle moves along the circle of the 
radius $\textfrac{1}{|p|}$. The distance from the starting point to the ending 
point is changing from  0 to $\textfrac{2}{|p|}$. 
For  $p$ large enough the geodesics of the length  $s$ 
belong to the interior of the cone with the angle of inclination 
 $\textfrac{|\varphi|}{4} s$ and the apex at the origin. 

\section{Distance along the  $x^4$ coordinate}

In the subspace $(x^2,x^3)$ the particle moves along the circle. 
Let us find the values of  $t$ for which the end of the geodesic is over the 
starting point 
in the subspace $(x^2,x^3,x^4)$.
The end of the geodesic is over the starting point for 
$pt=2\pi k$, $k\in\mathbb Z$.
Use the Taylor expansion for the equation
\eqref{x2x3x4} at the point 
$p=\textfrac{2\pi k}{t}+\theta$, $\theta\to 0$:
\begin{equation}
x^4(t)=-\varphi\frac{t^2}{4\pi k}
+\varphi\frac{t^3}{4\pi^2 k^2} \theta 
+O(\theta^2) \quad (\theta\to 0).
\end{equation}
Since $t=s$, the distance for which the particle moved along the  $x^4$ 
coordinate
is of the order  $\varphi\textfrac{s^2}{4\pi k}$.

For distribution with the positively defined metric tensor there is the theorem
about the ball shape \cite{Ger,Trelat}. In the coordinate neighbourhood of any 
point
one can choose the vector fields  $e_i$ which form the basis of the distribution 
$\mathcal A$. The field of planes that span commutators of all possible vector 
fields of the distribution  $\mathcal A$ is designated 
$\mathcal A_2=[\mathcal A,\mathcal A]$.
One also consider  $\mathcal A_{k+1}=[\mathcal A_k, \mathcal A]$.
This sequence is stabilizing: there is the  minimum $m$ such that 
$\mathcal A_m = \mathcal A_{m+1}$.
The number $m$ which can depend of the point of the manifold 
is called the degree of nonholonomity of the distribution $\mathcal A$.
Designate $n_k=\dim \mathcal A_k$ and define the function
$\varphi(i)=j$, if $n_{j-1}<i\leq n_j$.
The function $\varphi$ is defined on the set of numbers 
$\{1,\ldots\!,\dim \mathcal A_m\}$
with values at $\{1,\ldots\!,m \}$.
For completely nonholonomic distribution on the Riemannian manifold there are
the coordinates  $x^i$ such that
the ball of accessibility is bounded from above and from below by the set
$|x^i|\leq \varepsilon^{\varphi(i)}$.
For the two-dimensional distribution on the 3-dimensional manifold
$n_0=0$, $n_1=2$, $n_2=3$.
Therefore the function $\varphi(1)=1$, $\varphi(2)=1$,
$\varphi(3)=2$.

\section{Lagrange formulation for the equations of motion}
\subsection{General case} 

Let us consider the classical problem of the calculus of variations:
how to find an absolutely continuous vector-function which
maximizes the functional
\cite{VerGer,Pont}
\begin{equation}
J(x(\cdot),u(\cdot))=\int_0^T L(t,x(t),u(t))\, dt,
\end{equation}
where $u(t)=\textfrac{dx}{dt}$.
Assume that $x(0)=x_0$, $x(T)=x_1$, 
and allowed paths satisfy the conditions
$\varphi^i(t,x,\textfrac{dx}{dt})=0$, $i=1,\ldots\! k$.
Then there are 
$k$ measurable and limited functions $\lambda_i(t)$
called Lagrange multipliers and
a constant $a_0\geq 0$, at least one non-zero and such that
the function
\begin{equation}
L_\lambda(t,x(t),\dot x(t))=a_0 L(t,x(t),\dot x(t))+
\sum_{i=1}^k \lambda_i(t) \varphi^i(t,x(t),\dot x(t))
\end{equation}
almost anywhere on $[0,T]$ satisfies the Euler -- Lagrange equations
in the integral form
\begin{equation}
\dd{L_\lambda(t,x(t),\dot x(t))}{\dot x^i(t)}=
\int_0^T \dd{L_\lambda(\tau,x(\tau),\dot x(\tau))}{x^i}\, d\tau 
+ c_i, \quad i=1,\ldots\! n,
\end{equation}
where $c_i$ are constants \cite[p. 279]{Pont}.
The parameters $a_0$, $\lambda_i$ can be multiplied by any positive constant.
With $a_0=1$ we obtain regular geodesics.
With $a_0=0$ we obtain abnormal geodesics.

Now assume (for the general case) that the distribution $\mathcal B$ is defined
by the family of differential forms
$\omega^j$, $j=m{+}1,\ldots\! n$.
The horizontality conditions have the form 
\begin{equation} \label{horizonal}
\omega^j(u) = 0, 
\quad j=m{+}1,\ldots\! n, \quad \dot x = u.
\end{equation}
Then the Lagrange function for the considered problem is
\begin{equation}
L_{\lambda}(x,u) = a_0 \langle u, u \rangle^{1/2} + 
\sum_{j=m+1}^n \lambda_j \omega^j(u),
\end{equation}
where $\langle\cdot,\cdot\rangle$ is some non-degenerate\footnote{
This assumption for our model is not necessary, as shown below.
It is enough for the {\it restriction} of the metric tensor on the distribution
to be non-degenerate.} 
bilinear form on $TM^n$.
The equations of the horizontal geodesics are
\begin{equation}\label{EuLagr}
a_0 \langle \frac{Du}{dt}, \cdot\,\rangle +
\sum_{j=m+1}^n \dot\lambda_j \omega^j +
\sum_{j=m+1}^n \lambda_j d\omega^j(u,\cdot) = 0,
\end{equation}
where $\textfrac{Du}{dt}$ is the covariant derivative of the velosity vector
of the particle along the path.
The velocity vector can be decomposed by the basis vector fields of the distribution:
\begin{equation}
u(t)=\sum_{k=1}^m v_k(t) \xi_k(x(t)).
\end{equation}
Then the covariant derivative is
\begin{equation}\label{dudt}
\frac{Du}{dt}=
\sum_{k=1}^m \frac{dv_k}{dt} \xi_k + 
\sum_{i,j=1}^m v_i v_j \nabla_{\xi_i} \xi_j,
\end{equation}
where $\nabla$ is the symmetric connection defined by the bilinear form.
Use the Maurer -- Cartan formula
\begin{equation}
d\omega(\xi,\eta) = \xi\omega(\eta) - \eta\omega(\xi) -
\omega([\xi,\eta]).
\end{equation}
Since $\omega^j(\xi_i) = 0$, $i=1,\ldots\! m$, $j=m{+}1,\ldots\! n$,
then for the velocity vector and the basis of the distribution
$d\omega^j(u,\xi_i) = -\omega^j([u,\xi_i])$.
This commutator can be rewritten using the structural constants of the distribution
$\mathcal B$:
\begin{equation}
[u,\xi_i] = \sum_{k=1}^m v_k [\xi_k,\xi_i] =
\sum_{k=1}^m v_k \sum_{s=1}^n c^s_{ki} \xi_s,
\end{equation}
where $[\xi_k,\xi_i] = \sum\limits_{s=1}^n c^s_{ki} \xi_s$.
This basis and the differential forms can be selected such that
$\omega^j(\xi_i) = \delta_i^j$, $i,j=m{+}1,\ldots\! n$.
Then for the distribution
$\mathcal B$
\begin{equation}\label{eqnsA}
a_0 \langle \frac{Du}{dt}, \xi_l\rangle +
\sum_{j=m+1}^n \lambda_j \sum_{k=1}^m c^j_{lk} v_k = 0,
\quad l=1,\ldots\! m.
\end{equation}
Since 
$\omega^j(\xi_i) = \delta_i^j$, $i,j=m{+}1,\ldots\! n$,
then for the velocity vector and all basis vector fields of the tangent space
$d\omega^j(u,\xi_i) = -\omega^j([u,\xi_i])$.
Hence we can write the projection of the equation 
\eqref{EuLagr} at the basis vectors which do not belong to the distribution:
\begin{equation}\label{eqnsO}
a_0 \langle \frac{Du}{dt}, \xi_l\rangle +
\dot\lambda_l +
\sum_{j=m+1}^n \lambda_j \sum_{k=1}^m c^j_{lk} v_k = 0,
\quad l=m{+}1,\ldots\! n.
\end{equation}

\subsection{Lagrange formulation for the distribution $\mathcal A$
in our model}

The distribution $\mathcal A$ is defined by the differential form 
$\omega_x=\sum\limits_{i=0}^3 A_i(x) dx^i + dx^4$.
The Euler -- Lagrange equations for the length functional
\eqref{length} with the condition $\omega_{\gamma(t)}(\gamma'(t))=0$
are
\begin{equation}
a_0 \langle\frac{D\dot\gamma}{dt},\cdot\,\rangle +
\dot\lambda \omega + \lambda d\omega(\dot\gamma,\cdot\,) = 0,
\end{equation}
where the covariant derivative 
\begin{equation}
\frac{D\dot\gamma}{dt} =
\sum_{k=0}^3 \frac{dv^k}{dt} e_k + \sum_{i,j=0}^3 v^i v^j \, \nabla_{e_j} e_i,
\end{equation}
$\nabla$ is the symmetric connection defined by the metric tensor.
Please note that in the Lagrange method we have to use the metric tensor defined
on the whole tangent space $TM$ whereas in the maximum principle we used
the restriction of the metric tensor on the distribution only.
This is no problem because in the result we use the components of 
the metric tensor of the distribution only.
We can also use the connection (Christoffel symbols) with all lower indexes
and in this case the bilinear form does not have to be non-degenerate.
The velocity vector 
$\dot\gamma(t)=\sum\limits_{k=0}^3 v^k(t) e_k(\gamma(t))$.
The Euler -- Lagrange method also has the constant
$a_0\geq 0$ \cite[p.~279]{Pont}. 
The parameters $a_0$, $\lambda$ can be multiplied by any positive constant.
With $a_0=1$ we obtain the equations of regular geodesics.
With $a_0=0$ we obtain the equations of abnormal geodesics.
Projecting the equations of motion on the basis of the distribution
$\mathcal A$ we obtain
\begin{equation}
a_0 \langle\frac{D\dot\gamma}{dt}, e_i \rangle +
\lambda \sum_{j=0}^3 c^4_{ij} v^j = 0,
\quad i=0,\ldots\!,3.
\end{equation}
The structural constants are defined by
$[\xi_i,\xi_j]=\sum\limits_{l=0}^4 c^l_{ij} \xi_l$, 
where $\xi_j=e_j$, $j=0,\ldots\!,3$, $\xi_4=\partial_4$.
For the distribution $\mathcal A$ $[e_i,e_j]=-F_{ij}\partial_4$,
$[e_j,\partial_4]=0$. Therefore $c^4_{ij}=F_{ji}$, and other structural 
constants are zero.
Since $\omega_x(\partial_4)=1$, then
\begin{equation}
a_0 \langle\frac{D\dot\gamma}{dt}, \partial_4 \rangle + \dot\lambda = 0.
\end{equation}
We can assume that $\langle\cdot,\partial_4\rangle=0$.
Then $\lambda=\text{const}$, and $\lambda$ can be interpreted as the charge of the particle
(or charge to mass ratio depending on the Lagrangian).
Note that applying the maximum principle in section \ref{proof} we used
the restriction of the metric tensor on the distribution only.
Since this restriction is non-degenerate, we can rise appropriate indexes
of the Christoffel sybols and get ordinary differential equations.
The result does not depend on the extension of the metric tensor
on the whole tangent space.

All figures published in this paper were produced using our own 3D graphics program
\copyright V.R.~Krym.


\begin{thebibliography}{29}

\bibitem{Gromov} Gromov M. {\it Carnot-Caratheodory Spaces
Seen From Within.} Preprint IHES/M/94/6 (Institut des Hautes
Etudes Scientifiques, 1994)

\bibitem{VerGer} Vershik A.M., Gershkovich V.Ya. 
{\it The Nonholonomic Dynamical Systems.
Geometry of Distributions and Variational Problems.}
Dynamical Systems--7. Results of Science and Technique, ser. "The 
Contemporary Problems of Mathematics, Fundamental Directions",
v. 16, pp. 5--85. Moscow, 1987. (Russian)

\bibitem{Dobr} Dobronravov V.V. Foundations of the Mechanics of
the Nonholonomic Systems. Moscow, 1970. (Russian)

\bibitem{Newm} Newmark Yu.I., Fufaev N.A. The Dynamics of the Nonholonomic
Systems. Moscow, 1967. (Russian)

\bibitem{Beem} Beem J., Ehrlich P. 
Global Lorentzian Geometry. 
Marcel Dekker, 1981.

\bibitem{Kr} Krym V.R. 
{\it Geodesics Equations for a Charged Particle in the Unified Theory
of Gravitational and Electromagnetic Interactions.}
// Teor.  Matem. Fisika, 1999, v. 119, N 3, pp. 517--528. (Russian)
// Theoretical and Mathematical Physics, 1999, 119:3, 811--820 
(English)


\bibitem{Landau} Landau L.D., Lifshitz E.M. Theoretical Physics, v. 2.
The Field Theory. Moscow, 1988. (Russian)


\bibitem{Gray} Gray C.G., Karl G., Novikov V.A. {\it
Progress in Classical and Quantum Variational Principles.}
Reports on Progress in Physics, 2004, v. 67, N2, pp. 159--208.

\bibitem{Rumer} Rumer Yu.B. Researches in 5-Optics.
Moscow, 1956. (Russian)

\bibitem{BailinL} Bailin D., Love A. {\it Kaluza -- Klein Theories.}
Reports on Progress in Physics, 1987, v. 50, pp. 1087--1170.

\bibitem{Sriv} Srivastava S.K. {\it Some Aspects of Kaluza--Klein Cosmology.}
Pramana--Journal of Physics, 1997, v. 49, N4, pp. 323--370.

\bibitem{KrM} Krym V.R. {\it Smooth Manifolds of Kinematic Type.}
// Teor.  Matem. Fisika, 1999, v. 119, N 2, pp. 264--281. 
(Russian)
// Theoretical and Mathematical Physics, 1999, 119:2, 605--617 (English)

\bibitem{KrP} Krym V.R, Petrov N.N. 
{\it Causal Structures on Smooth Manifolds.}
Vestnik Sankt-Peterburgskogo Universiteta, 
ser. 1, 2001, N 2, pp. 27--34. (Russian)

\bibitem{KrP2} Krym V.R, Petrov N.N. 
{\it Local Ordering on Smooth Manifolds.}
Vestnik Sankt-Peterburgskogo Universiteta, 
ser. 1, 2001, N 3, pp. 32--36. (Russian)

\bibitem{Vasil} Vasiljev F.P. 
Numerical Methods of Solution of Extremal Problems.
Moscow, 1988. (Russian)

\bibitem{KrotG} Krotov V.F., Gurman V.I.
Methods and Problems of Optimal Control. 
Moscow, 1973. (Russian)

\bibitem{Pont} Ponrjagin L.S., Boltjanskij V.G., 
Gamkrelidze R.V., Mishchenko E.F.
The Mathematical Theory of Optimal Processes.
Moscow, 1983. (Russian)

\bibitem{Phil} Filippov A.F. 
{\it On Some Questions of the Theory of Optimal Control.}
Vestnik Moskovskogo Universiteta, ser. Math. \& Mech.,  1959, N 2, pp. 25--32. (Russian)

\bibitem{Petr} Petrov N.N. 
{\it Existence of Abnormal Shortest Paths in sub-Riemannian Geometry.}
Vestnik Sankt-Peterburgskogo Universiteta, 
ser. 1, 1993, Iss. 3 (N 15), pp. 28--32. (Russian)

\bibitem{Bonn} Bonnard B., Chyba M.
{\it Singular Trajectories and Their Role in Control Theory.}
Mathematiques \& Applications, v. 40. Paris: Springer, 2003.

\bibitem{Mont1} Montgomery R. {\it A Survey of Singular Curves in
sub-Riemannian Geometry.}
J. Dynam. Contr. Syst., 1995, v. 1, pp. 49--90.

\bibitem{Mont2} Montgomery R. {\it Survey of Singular Geodesics.}
Progress in Math., 1996, v. 144, pp. 325--339.

\bibitem{Kup} Kupka I. 
{\it Sub-Riemannian Geometry.}
Asterisque, 1997, v. 241, pp. 351--380.

\bibitem{KupO} Kupka I., Oliva W.M.
{\it The non-Holonomic Mechanics.}
J. Differ. Equ., 2001, v. 169, N 1, pp. 169--189.

\bibitem{Ger} Gershkovich V.Ya. 
{\it Two-side Estimations of Metrics Generated by Absolutely Nonholonomic
Distributions on Riemannian Manifolds.}
Doklady Akademii Nauk, 1984, v. 278, N 5, pp. 1040--1044. (Russian)

\bibitem{Trelat} Trelat E. 
{\it Non-Subanalyticity of sub-Riemannian Martinet Spheres.}
C. R. Acad. Sci. Paris, Ser. I, Math., 2001, v. 332, N 6, pp. 527--532.

\bibitem{Mitch} Mitchell J. {\it On Carnot-Caratheodory Metrics.}
J. Diff. Geom., 1985, v. 21, N 1, pp. 35--45.

\bibitem{Jean} Jean F. {\it Uniform Estimation of sub-Riemannian Balls.}
J. Dyn. Control Syst., 2001, v. 7, N 4, pp. 473--500.

\bibitem{Griff} Griffiths Ph. A. Exterior Differential Systems and the 
Calculus of Variations. Birkhauser, 1983.
(Progress in Mathematics, v. 25).

\end{thebibliography}
\end{document}